\documentclass[12pt]{article}
\usepackage{a4wide}
\usepackage{amssymb}
\usepackage{graphicx}
\newtheorem{thm}{Theorem}
\begin{document}
\title{Visualizing Ricci Flow of Manifolds of Revolution}
\author{J. Hyam Rubinstein
and
Robert Sinclair \\
Department of Mathematics \& Statistics \\
The University of Melbourne \\
Parkville \\
Victoria 3010 \\
Australia}
\maketitle

\centerline{\bf Abstract}
{\small We present numerical visualizations of Ricci Flow of surfaces
and $3$-dimensional
manifolds of revolution.
{\tt Ricci\_rot} is an educational tool which visualizes surfaces of
revolution
moving under Ricci
flow. That these surfaces tend to remain embedded in $\mathbb{R}^3$
is what makes
direct visualization possible. The numerical lessons gained in
developing this tool may
be applicable to numerical simulation of Ricci flow of other
surfaces. Similarly for simple
3-dimensional manifolds like the $3$-sphere, with a metric which is
invariant under the action of $SO(3)$ with
$2$-sphere orbits, the metric can be represented by a $2$-sphere of
revolution, where the distance to the axis of revolution represents
the radius of a $2$-sphere orbit. Hence we can also visualize the
behaviour of such a metric under Ricci flow. We discuss
briefly why surfaces
and $3$-manifolds of revolution remain embedded in $\mathbb{R}^3$ and
$\mathbb{R}^4$
respectively, under Ricci flow and finally
indulge in some speculation
about the idea
of Ricci flow
in the larger space of positive definite and indefinite metrics.}

\section{Introduction}

Understanding Ricci flow has become fundamental to attempts to prove
Thurston's geometrization conjecture (\cite{CAO}, \cite{CHOW},
\cite{PEREL1}, \cite{PEREL2}). Since Ricci flow is geometric by
nature, it is
natural to ask whether it can be visualized.
Unlike mean-curvature and related flows, which act in the first
instance as a force pushing
on a surface
and therefore tend to keep surfaces embedded in Euclidean space, Ricci
flow
acts directly on the metric of the surface, tending not to preserve
embeddedness.
A number of interesting results have been obtained by restricting to
classes of metrics of revolution, since such symmetries
are preserved under Ricci flow and the metric depends on considerably
fewer parameters in such cases (\cite{SIMON}, \cite{ANG},
\cite{IVEY}). We will discuss for this class of metrics,
why Ricci flow does preserve embeddedness in Euclidean space (\cite{ENGMAN}).
Moreover, the key phenomenon of neck pinching occurs naturally
for metrics of revolution. Good estimates are available to understand
the limiting behaviour of metrics as they become singular under neck
pinching in (\cite{SIMON}, \cite{ANG}). Our visualization approach
enables
us to generate interesting pictures and information about neck pinching.
In \cite{BCHOW}, it is shown that any metric on the $2$-sphere
converges to the round metric
under Ricci flow, so neck pinching occurs only in higher dimensions.

Note that there are many interesting questions about Ricci flow for
$3$-dimensional manifolds. There are a limited number of examples
where one knows
completely how the metric evolves (\cite{CHOW}). It would be very
useful to have further examples where Ricci flow continues for
infinite time, without developing
singularities. For example, there are important constructions (Dehn
surgery) of metrics with strictly negative sectional curvature on
$3$-manifolds (\cite{GT},
\cite {BH}). But even in this case, it is not clear whether
singularities will form under Ricci flow. In fact, no examples are
known where a metric with
all sectional curvatures negative flows to a metric with positive
sectional curvature at some points.

\section{Preserving the Embedding into Euclidean Space}

In this section, we give a quick discussion of why $n$-dimensional
spheres of revolution ($K_2=1$ -- see Section \ref{HASS}), which are
initially isometrically
embedded into Euclidean $(n+1)$-dimensional space, remain isometrically
embedded, so long as Ricci flow has a smooth solution.

\subsection{Characterising when Isometric Embeddings exist}

We observe a simple extension of \cite{ENGMAN}. Note that in later
sections, we will use extrinsic coordinates, which are
useful to see how points of the surface or manifold move under Ricci
flow when visualized in Euclidean space. Here it is more
convenient to use intrinsic coordinates. So on $S^n$, a metric
is chosen of the form  $ds^2+\alpha^2(s)d\theta^2$,
where $s$ is arc length along a geodesic joining the North Pole to
the South Pole and $\theta$ represents the coordinates in a
hypersphere of radius $\alpha(s)$, which is the orbit of an isometric
action of $SO(n)$ on $S^n$, where the orbit space is an arc.
We assume that $s$ varies between $0$ and $L$.

In \cite{ENGMAN}, it is shown that a necessary and sufficient
condition that $S^2$, with a metric of this form, can be isometrically
embedded in $\mathbb{R}^3$ as a surface of revolution, is that the
integral of Gaussian curvature over any ``polar cap'' is positive,
where such
a region consists of all points
where either $0\le s\le k$ or $k \le s \le L$, for some $k$ with
$0<k<L$. Note also that a surface of revolution in $\mathbb{R}^3$ is
invariant under
an isometric action of $SO(2)$ by rotation about the $x^1$ axis.
Similarly we define a manifold of revolution in $\mathbb{R}^{(n+1)}$
as being invariant under the isometric action of $SO(n)$ by rotation
about the $x^1$ axis. Finally we follow \cite{ENGMAN}, by requiring that
any manifold of revolution $\mathbb{R}^{(n+1)}$ has at most one component of
intersection with any hyperplane of the form $x^1=k$ for any constant $k$.
Equivalently, the manifold is obtained by rotating the graph of a
function of one variable
around the $x^1$ axis, by the action of $SO(n)$.
Note that for visualization purposes, it is not essential but is
certainly convenient to have this
restriction.

\begin{thm}
Suppose that a metric of the form $ds^2+\alpha^2(s)d\theta^2$ is
chosen on $S^n$. Then there is an isometric embedding as a manifold
of revolution in $\mathbb{R}^{(n+1)}$ if and only if the same
metric, viewed as on $S^2$, can be isometrically embedded into
$\mathbb{R}^3$ as
a surface of revolution.

\end{thm}

The proof of this theorem is very easy and is left to the reader.

\subsection{Ricci Flow preserves Embeddability of Manifolds
of Revolution}

An observation about the characterisation in \cite{ENGMAN} is
as follows. For a polar cap, Gauss Bonnet shows that the sum of the
integral of the Gaussian curvature of the cap and the integral of the curvature
of the boundary curve $C$ is $\pi$. Therefore positivity of the first
integral is equivalent to the second integral being strictly less
than $\pi$.
Note that $C$ is the orbit of $SO(2)$ acting on a plane parallel to
the $x^2 x^3$ plane and so is a standard round circle in this plane,
with centre on
the intersection of the plane with the $x^1$ axis. Hence the integral
of its curvature in $\mathbb{R}^{3}$ is exactly $\pi$ and the
direction of its curvature at
each point points towards the $x^1$ axis. But then to compute the
curvature in the surface of revolution, one has to project onto the
tangent space of this surface.
Hence the projected curvature vector will be shorter, unless the
tangent plane is vertical at every point of $C$. So this shows that
the condition on the polar caps remains true as
Ricci flow proceeds, unless at some time and for some value of $s$,
with $0<s<L$,
$|\alpha^{\prime}(s)|=1$. Therefore we need to show that
$|\alpha^{\prime}(s)|<1$ remains true, away from the poles, so long
as Ricci flow
produces a smooth solution (\cite{ANG},\cite{ENGMAN},\cite{IVEY}),
to prove the following result.

\begin{thm}
If there is an isometric embedding of $S^n$ as a manifold
of revolution in $\mathbb{R}^{(n+1)}$ then the manifold remains
isometrically embedded so long as Ricci flow gives a smooth solution.

\end{thm}

{\bf Proof}
Let $v= \alpha^{\prime}$, where $v$ and $\alpha$ are viewed as functions
of $s$ and time $t$, and as above, the derivative of $\alpha$ is taken
with respect to $s$.
By equation (16) of \cite{ANG}, $v$ satisfies the following evolution equation.
$${\partial_t}v = v_{ss} +{n-2 \over \alpha}v{v_s} + {{n-1} \over
{\alpha^2}}(1-v^2)v.$$
Here we are using $v_s$ to denote partial differentiation
in the $s$ variable.

Now initially, $|v|=|\alpha_s|<1$ for all values $0<s<L$,
by assumption. Let $t$ be a first value of time for which
$|v(s,t)|=1$, at some $s$ not corresponding to the North or South
pole of the sphere. If $v(s,t)=1$, then this is a maximum value
of $\alpha$ whereas if $v(s,t)= -1$, then this is a minimum value.

So we see by
the maximum principle, (see eg \cite{CHOW}) that this
gives a contradiction, since $1-v^2=v_s=0$ for both cases. This shows
that the sphere does remain embedded,
so long as it is smooth.

\section{Computational Formulation}

A two-dimensional surface of revolution of genus zero embedded in
$\mathbb{R}^3$ can be
defined in a polar representation with coordinates $x^1=\rho$ and
$x^2=\theta\in[0,2\pi[$ by a metric of the form
\[
[g_{\mu\nu}]=\left[\begin{array}{cc}
h(\rho) & 0 \\
0 & m(\rho) \end{array}\right],
\]
where $\sqrt{m(\rho)}$ has the direct physical interpretation as the
radius from the axis
of rotation. For a closed surface, we require
\[
m(\rho_{\scriptstyle pole})=0.
\]
We have chosen $\rho_{\mbox{\footnotesize North\_Pole}}=0$
and $\rho_{\mbox{\footnotesize South\_Pole}}=\pi$.
When $h(\rho)$ is a constant, $\sqrt{h}\times\rho$ is the distance of
a point with
coordinates $(\rho,\theta)$ from the North Pole along a meridian.
Smoothness at the poles demands
\begin{equation}
\left.
\frac{\partial\sqrt{m(\rho)}}
{\partial\rho}\right|_{\rho=\rho_{\scriptscriptstyle
pole}}
=\sqrt{h(\rho_{\scriptstyle pole})},
\qquad\rho_{\scriptstyle pole}\in\{0,\pi\}.
\label{SMOOTH_ENDS}
\end{equation}
Note that the necessary and sufficient condition for isometric
embeddability as a classical
surface of revolution in $\mathbb{E}^3$ for $S^1$-invariant metrics on
$S^2$ --
the nonnegativity of the integrals of the curvature over all disks
centred at
a pole --
is given in \cite{ENGMAN}.
The generating curve (cross-section)
of the surface of revolution is given by
\begin{eqnarray}
x(\rho) & = &
\int\limits_0^{\rho}\sqrt{h(s)-\left(\frac{\partial\sqrt{m(s)}}
{\partial
s}\right)^2}\,ds
\nonumber \\
y(\rho) & = & \sqrt{m(\rho)}
\label{GEN_CURVE}
\end{eqnarray}
with $\rho\in[0,\pi]$.

The Ricci tensor is
\begin{equation}
[R_{\mu\nu}]=\left\{\begin{array}{cl}
\left[\begin{array}{cc}
\displaystyle \frac{(m')^2}{4m^2}-
\frac{m''}{2m}+
\frac{m'\,h'}{4m\,h} & \displaystyle 0 \\[2ex]
\displaystyle 0 &
\displaystyle \frac{(m')^2}{4m\,h}-
\frac{m''}{2h}+
\frac{m'\,h'}{4h^2}
\end{array}\right] & (\,\rho\in\,]0,\pi[\,) \\[7ex]
\left[\begin{array}{cc}
\displaystyle
\frac{m''''}{4m''}-
\frac{h''}{2h} & \displaystyle 0 \\[2ex]
\displaystyle 0 &
\displaystyle 0
\end{array}\right] & (\,\rho\in\{0,\pi\}\,),\end{array}\right.
\label{RICKY_1}
\end{equation}
where primes denote partial differentiation with respect to $\rho$.

Ricci flow satisfies
\[
\frac{\partial g_{\mu\nu}}{\partial t}=-2R_{\mu\nu}.
\]

The possible initial metrics {\tt Ricci\_rot} provides are given by
\begin{equation}
\left[\begin{array}{cc}
\displaystyle \qquad1\qquad & \displaystyle 0 \\[2ex]
\displaystyle 0 &
\displaystyle \left(\frac{\sin
\rho+c_3\sin3\rho\,+c_5\sin5\rho}{1+3c_3+5c_5}\right)^2
\end{array}\right].
\label{C3C5}
\end{equation}
These surfaces are
symmetric under reflections about a plane normal to their axis of
rotation which
passes through their centre of mass. This extra symmetry has proven
useful in making
the simulation faster, since it then suffices to compute on only half
of the surface.

Attempting to integrate using Equation
(\ref{RICKY_1}) and central finite differences in a
purely explicit scheme leads, not
surprisingly to numerical instability, particularly (but not only) at
the poles.
This instability can be reduced significantly by noting that
\[
\frac{(m')^2}{4m^2}-\frac{m''}{2m}=\frac{-(\sqrt{m})''}{\sqrt{m}}.
\]

\subsection{Filtering and Reparametrization}

The numerical instabilities are however persistent enough to require
more potent
measures. We did not try any implicit scheme, which may have alleviated
the numerical instabilities, but would still have required very small
time steps near the poles.

Inspired by spectral methods \cite{SPECTRAL},
we first introduced a filter which consists of transforming to
Fourier space (DFT), dropping shorter
wavelength terms, and then transforming back. This process is
facilitated by noting that
the rotational and reflection symmetries of the family of surfaces
(\ref{C3C5}) are preserved
by Ricci flow.
We can therefore write
\[
h(\rho)=\sum\limits_{i=0}^{N_h}h_i\cos 2i\rho
\]
and
\[
\sqrt{m(\rho)}=\sum\limits_{i=0}^{N_m}m_i\sin\left((2i+1)\rho\right),
\]
where equality of course only applies to $N_h=N_m=\infty$.
We actually choose $N_h$ and $N_m$ to be much less than the number of
points used in the
finite difference scheme, and find that there is no need to implement
the FFT.
We need to ensure that (\ref{SMOOTH_ENDS}) is satisfied at the poles,
to be sure that the
generating curve (\ref{GEN_CURVE}) can be found. We do this by
multiplying the computed values
of $\sqrt{m}$ by
\[
\frac{\displaystyle
\frac{\sqrt{h(\rho_{\scriptstyle pole})}}
{\sum_{i=0}^{N_m}(2i+1)\times m_i}
+K\,(\rho-\rho_{\scriptstyle pole})^2
}{1+K\,(\rho-\rho_{\scriptstyle pole})^2},
\]
where $K$ is a small positive constant, but large enough that this
correction factor quickly
becomes unity away from the poles. Setting $K=0$ has the unwanted
side-effect that local
numerical errors near the poles can lead to global changes in shape
which do not satisfy the
Ricci flow.

By its very nature, the Ricci flow forces some parts of a surface to
contract while others
expand. This means that whatever set of nodes we have chosen for our
finite-difference scheme
will effectively become less and less evenly distributed, creating
further sources of
numerical instability where they become too dense. The solution to
this problem is to recall that
any given surface corresponds to many different parametrizations,
each with a corresponding
metric. From a numerical point of view, the most pleasant
parametrizations let
$h(\rho)$ be
a constant (this corresponds to $\rho$ being directly proportional to
distance from a pole).
Let
\[
\ell(\rho)=\int\limits_{s=0}^\rho \sqrt{h(s)}\,ds
\]
be the distance from the North Pole.
We have the freedom to reparametrize:
\[
[g_{\mu\nu}]=\left[\begin{array}{cc}
h(\rho) & 0 \\
0 & m(\rho) \end{array}\right]\mapsto
\left[\begin{array}{cc}
\displaystyle \left(\frac{\ell(\pi)}{\pi}\right)^2 & \displaystyle 0 \\
\displaystyle 0 & \displaystyle m\left(
\ell^{-1}\left(\frac{\rho\times\ell(\pi)}{\pi}\right)
\right) \end{array}\right]
\]
where $\rho\in[0,\pi]$ both before and after. It is this
reparametrization which,
performed after every few iterations of the explicit integration and
filtering steps, has made
{\tt Ricci\_rot} stable enough to use.

\section{Use of the Software}

First, a quick introduction:
Press ``{\tt n}'' and then drag the mouse to
choose an initial surface. Then hold down the up-arrow key to watch the
flow
until it stops. Once one has begun flowing, dragging the mouse will
rotate the
surface.

\subsection{Details}

{\tt Ricci\_rot} is an ANSI C program \cite{KERNI} which uses the
OpenGL standard
\cite{OPEN}
to display graphics.
It should therefore be portable, although it has only been tested on
a Mac OS X platform.

Upon launch, a window is opened with the image of a sphere with blue
meridians and parallels.
This initial shape can be deformed by holding down the left mouse
button and dragging in any
direction. The possible initial shapes belong to the family (\ref{C3C5}),
where $c_3$ is varied by horizontal, and $c_5$ by vertical motion.
The program will not accept all possible values of $c_3$ and $c_5$,
so, dragging
very far in any direction, one may find that the shape ceases changing.

To examine the shape, pressing ``{\tt f}'' (for ``flow'') will put the
program into flow
mode. The meridians and parallels change colour. Now, dragging the mouse
rotates the surface. To make any more changes, pressing ``{\tt n}'' (for
``new shape'') will put the program back in its initial mode. The
meridians and parallels change colour back to blue, and dragging the
mouse
changes the shape.

When one has chosen a shape one wishes to flow, it is enough to begin pressing
the up-arrow key. Usually this key will repeat if held down. Then
one can see
the surface flow continuously. Since the program is now in flow mode, dragging
the mouse will rotate the surface.

Once the program has met with a numerical instability, the flow
will stop, and
the parallels and meridians will appear black. Any ripples which may
appear
on the surface at this stage are a result of the numerical instability.

At any time, pressing the left-arrow or right-arrow keys will rotate
the surface
anti-clockwise or clockwise respectively.

The down-arrow key flows the surface backwards in time. This evolution is
highly unstable. One should not expect to be able to flow for very long.

Pressing ``{\tt m}'' at any stage will change the display mode from
showing
the surface to showing the components of the metric ($g_{1\,1}=h$ in
green and $g_{2\,2}=m$ in blue) and the surface's cross-section in white.
To revert to displaying the surface in $\mathbb{E}^3$, press ``{\tt s}''.

There is a bright spotlight. This can be turned on or off by pressing the
right mouse button. A menu will appear, from which one can choose to turn
the spotlight on or off (or change to and from ``flow'' and ``new
shape'' modes).
It can be rotated (in a rather non-intuitive way) by dragging the mouse
with the middle mouse button held down.

\subsection{Examples}

\begin{figure}[p]
\begin{center}
\includegraphics[width=\textwidth]{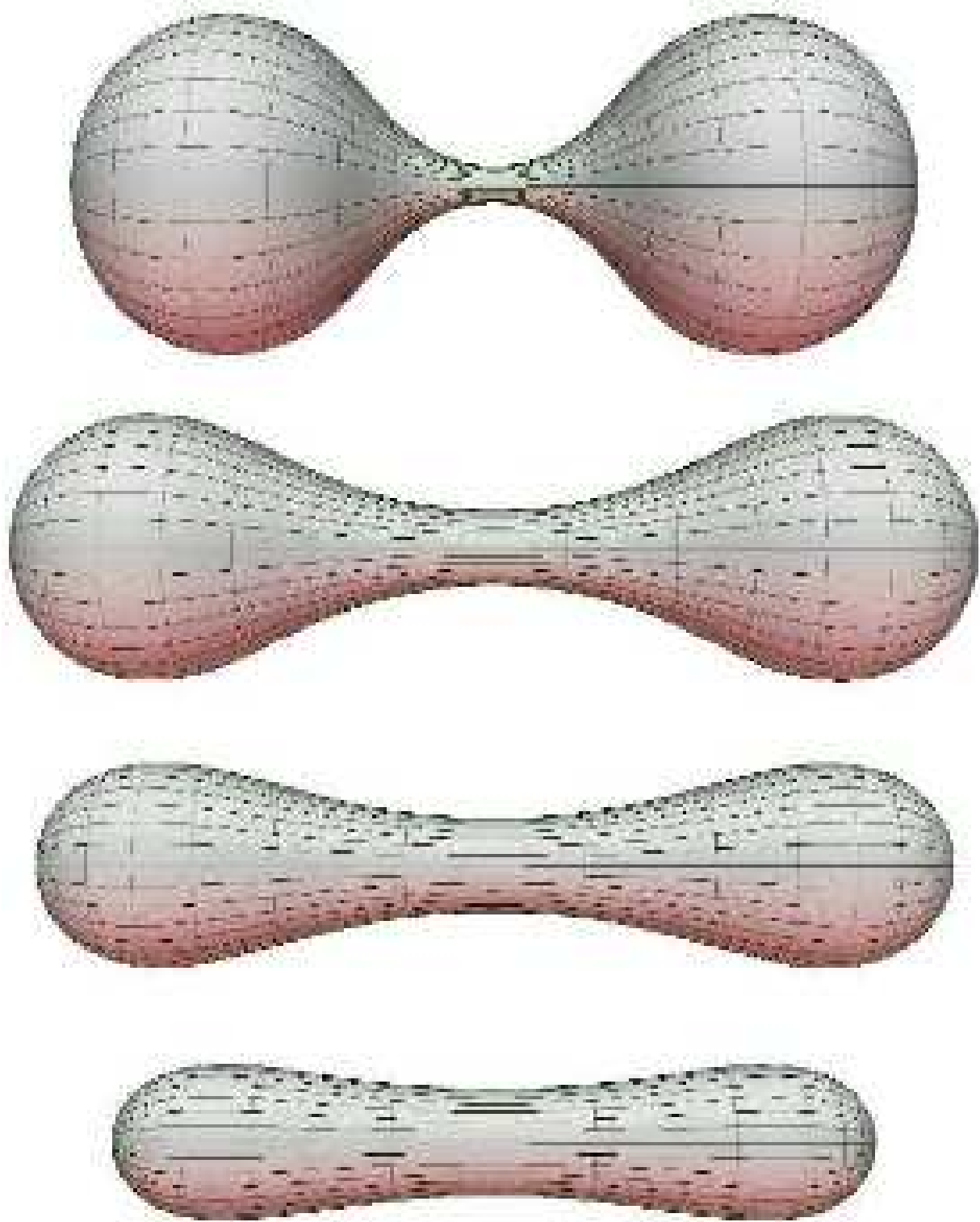}
\end{center}
\caption{\small Deformation of a dumbbell-like surface of revolution under
Ricci flow.
The pictures have been taken at equal time intervals and are drawn to the
same scale.}
\label{DUMB}
\end{figure}

\begin{figure}[tb]
\begin{center}
\includegraphics[width=0.95\textwidth]{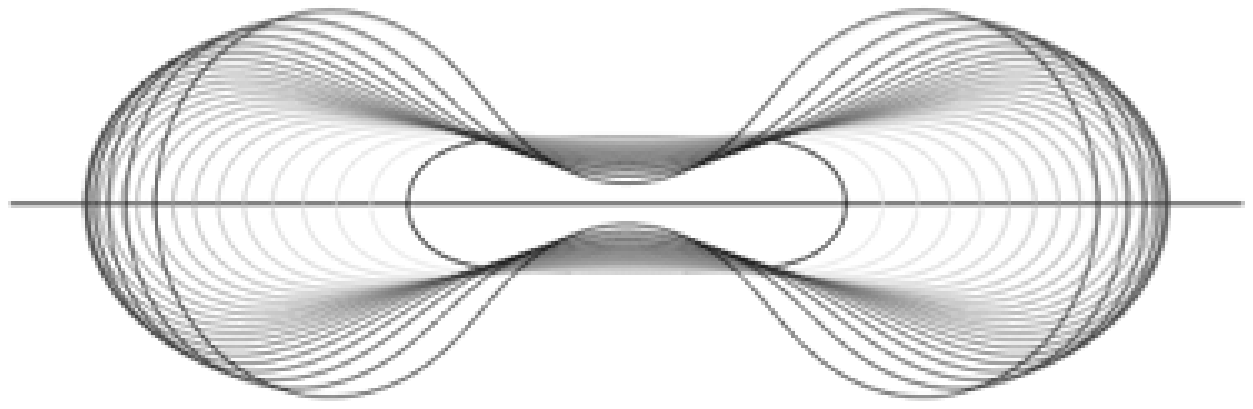}
\end{center}
\caption{\small Deformation of the cross-section of
a dumbbell surface of revolution under Ricci flow.
The axis of rotation is horizontal.
The pictures have been taken at equal time intervals. The initial and
final curves are both shaded black.}
\label{GENC_DUMB}
\end{figure}

The first example is of a dumbbell shape, given by Equation (\ref{C3C5}) with
$c_3=0.766$ and $c_5=-0.091$. This surface's flow is
illustrated in Figure \ref{DUMB} with time step $dt=0.01$.
The flow of its cross-section is
provided in Figure \ref{GENC_DUMB} time step $dt=0.002$.
It is interesting to compare this with Figure 4 of
\cite{MEAN}, which is an illustration of the mean-curvature flow of
the cross-section of what is initially a dumbbell shape.

\begin{figure}[p]
\begin{center}
\includegraphics[width=0.96\textwidth]{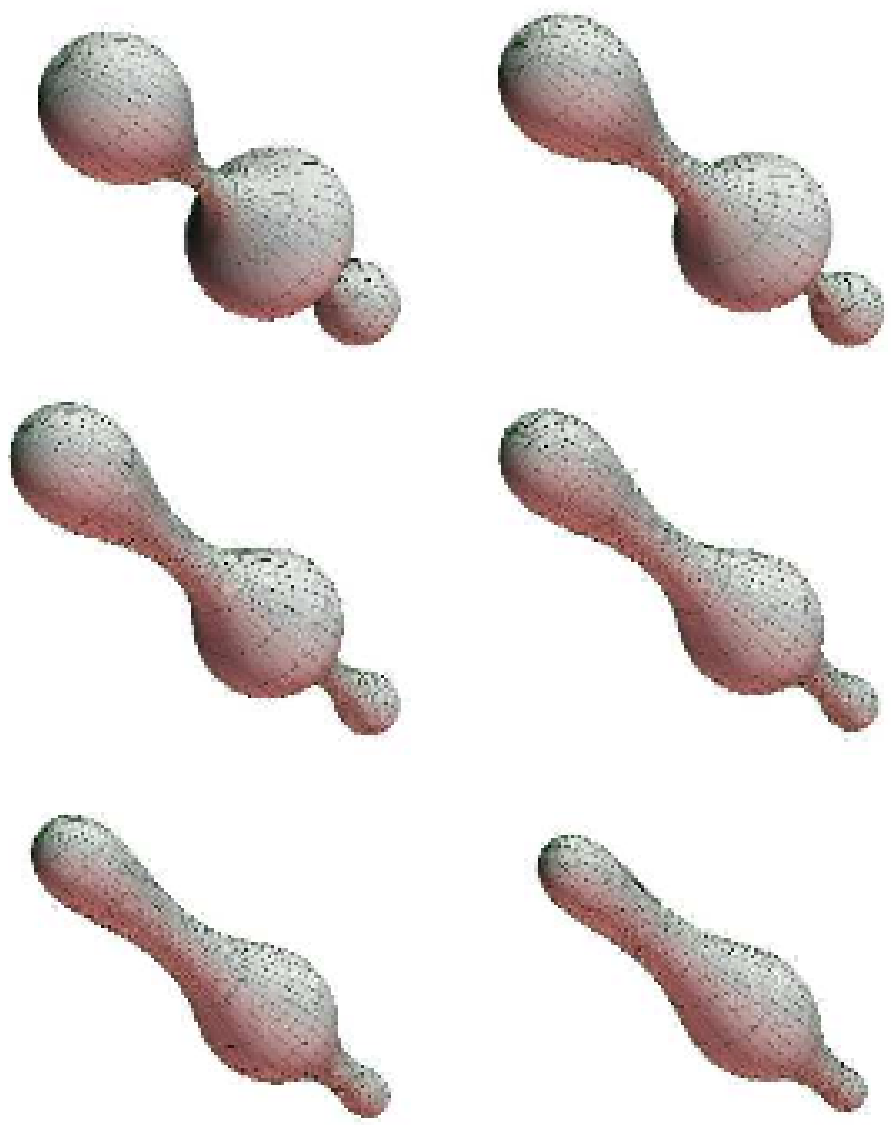}
\end{center}
\caption{\small Deformation of a surface of revolution under Ricci flow.
The pictures have been taken at equal time intervals and are drawn to the
same scale. The cross-sections are plotted in Figure \ref{GENC}.}
\label{FIG1}
\end{figure}

\begin{figure}[p]
\begin{center}
\includegraphics[width=0.92\textwidth]{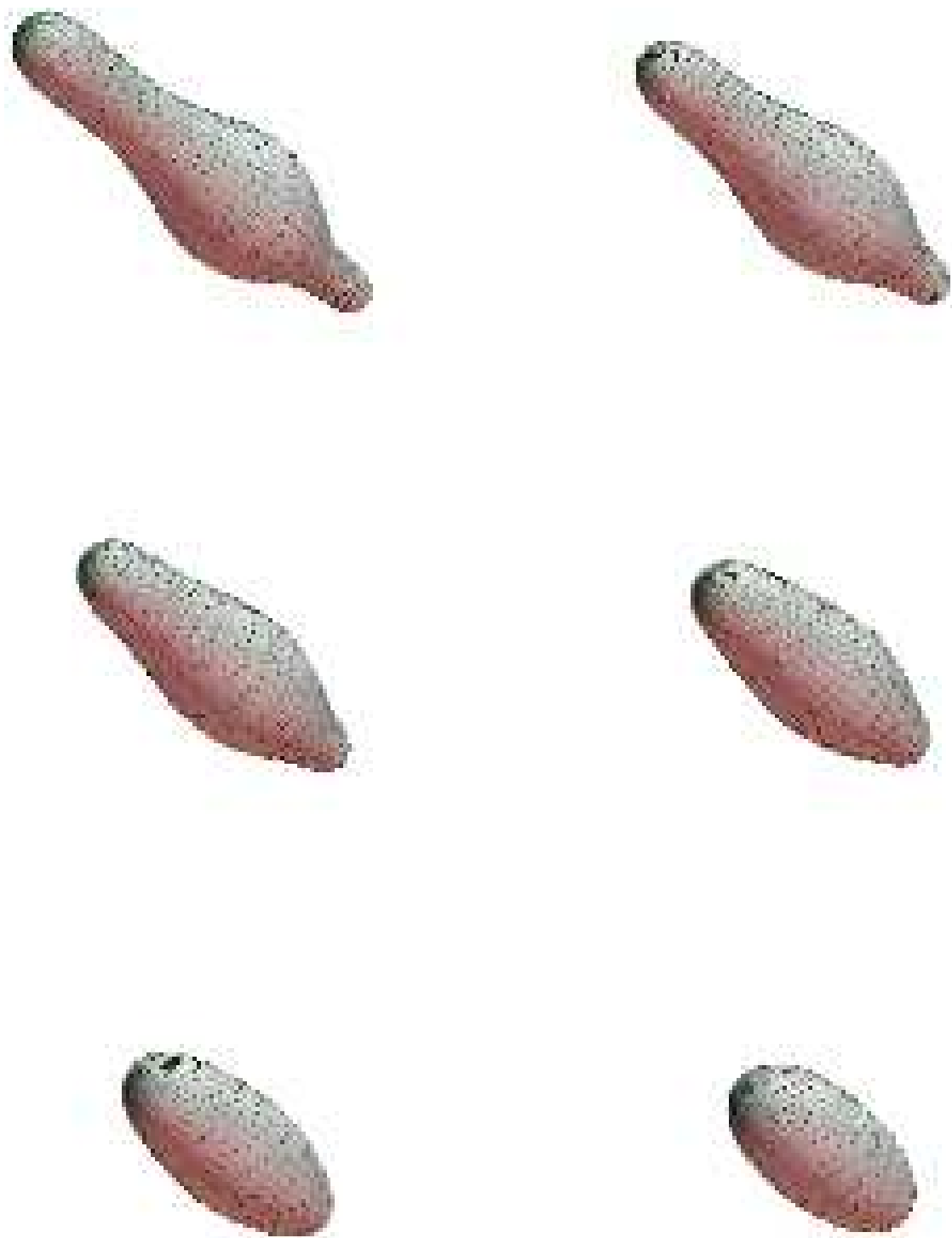}
\end{center}
\caption{\small Deformation of a surface of revolution under Ricci flow.
The pictures have been taken at equal time intervals and are drawn to the
same scale. These pictures are the continuation of the evolution of the
surface begun in Figure \ref{FIG1}.
The cross-sections are plotted in Figure \ref{GENC}.}
\label{FIG2}
\end{figure}

\begin{figure}[tb]
\begin{center}
\includegraphics[width=0.95\textwidth]{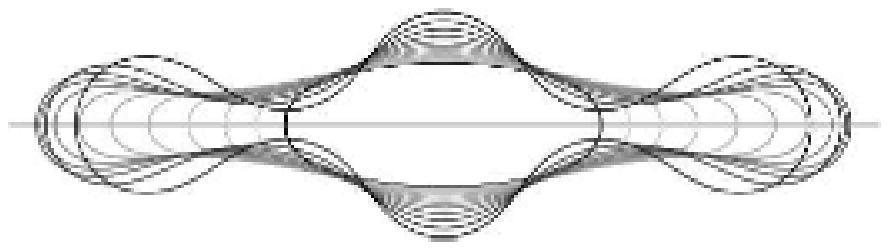}
\end{center}
\caption{\small Deformation of the cross-section of
a surface of revolution under Ricci flow.
The axis of rotation is horizontal.
The pictures have been taken at equal time intervals. The initial and
final curves are both shaded black.}
\label{GENC}
\end{figure}

\begin{figure}[t]
\begin{center}
\includegraphics[width=0.65\textwidth]{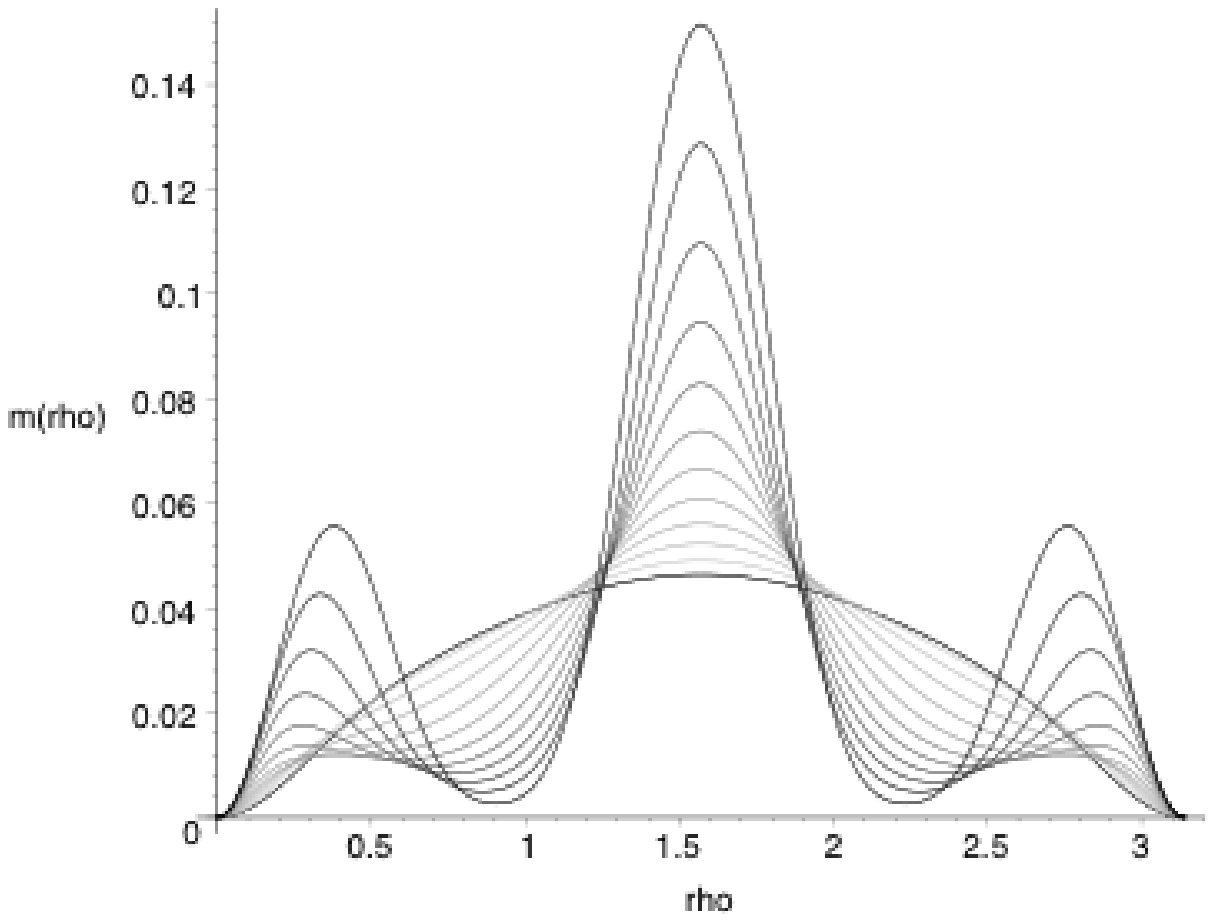}
\end{center}
\caption{\small Evolution of the metric component $g_{2\,2}=m(\rho)$ of
a surface of revolution under the combined action of Ricci flow
and reparametrization -- meaning that the corresponding curves
$g_{1\,1}=h$ are independent of $\rho$.
The curves have been plotted at equal time intervals. The initial and
final curves are both shaded black.}
\label{MC}
\end{figure}

Figures \ref{FIG1}, \ref{FIG2}, \ref{GENC} and \ref{MC} were
all generated using
the initial surface given by
$c_3=0.021$ and $c_5=0.598$, and with
the time step $dt=0.002$.
The values of $g_{1\,1}=h$ ($h(\rho)$ becomes a constant after
reparametrization)
are $1$, $1.027938$, $1.000526$, $0.936608$, $0.843907$,
$0.729106$, $0.601463$, $0.475008$, $0.365054$, $0.278860$,
$0.213761$ and $0.163754$.

The behaviour is in both cases what one expects -- flow towards spheres
of
constant positive Gaussian curvature \cite{BCHOW}. See also Chapter 5 of
\cite{CHOW}, in which it is shown that any solution of the unnormalized Ricci
flow on a topological $S^2$ shrinks to a round point in finite time,
following \cite{BCHOW}.
Note that the ``modified'' Ricci flow in \cite{HYAM}, for the special
case of surfaces
of revolution embedded in Euclidean space, could be visualized in a
similar way.

\section{Three-Dimensional Manifolds of Revolution}

\begin{figure}[t]
\begin{center}
\includegraphics[width=0.65\textwidth]{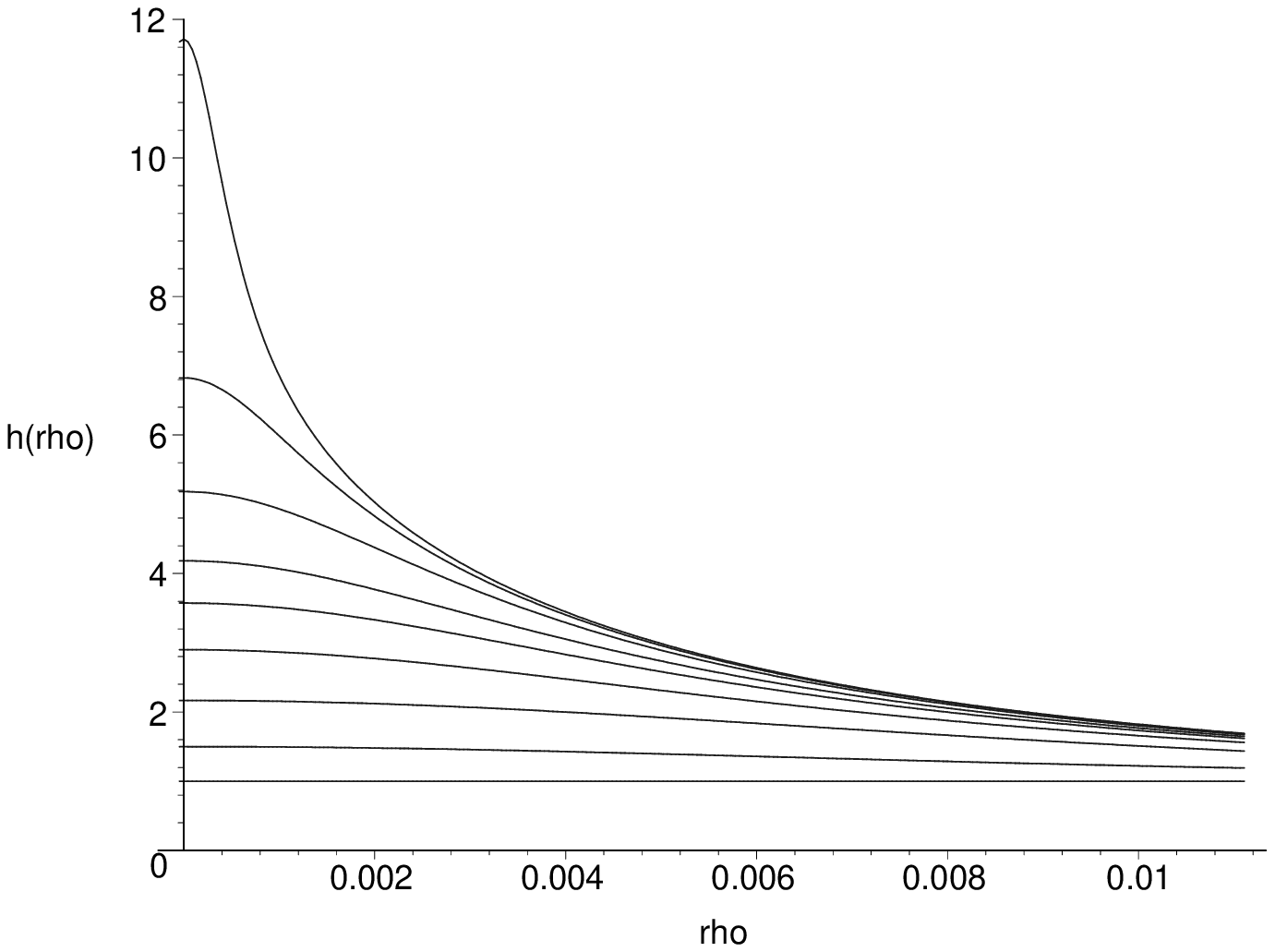}
\end{center}
\caption{\small Qualitatively correct
evolution of the metric component $g_{1\,1}=h(\rho)$ of
a 3-manifold of revolution under Ricci flow.
The curves have been plotted at unequal time intervals.}
\label{H_PLOT}
\end{figure}

\begin{figure}[t]
\begin{center}
\includegraphics[width=0.65\textwidth]{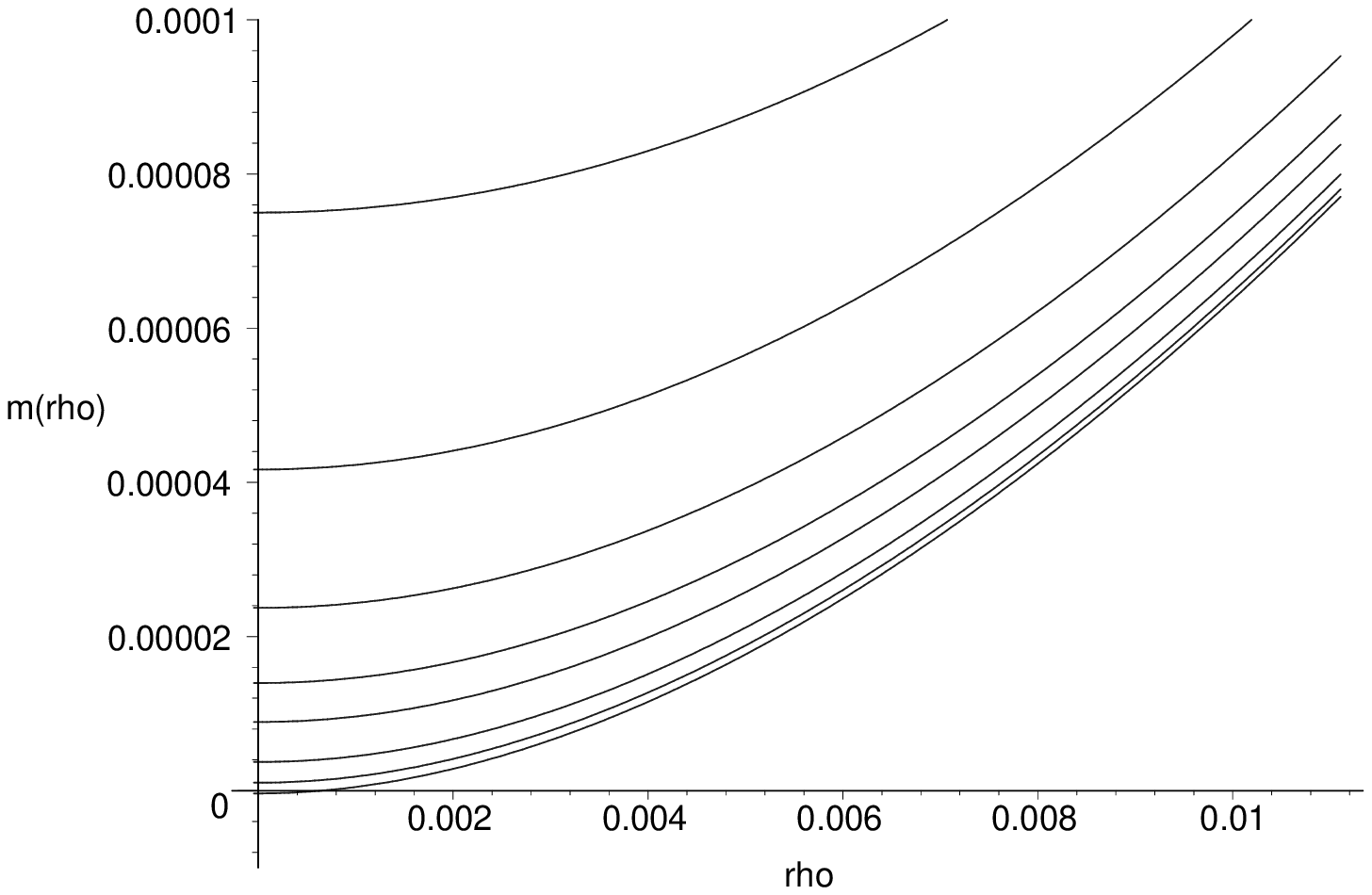}
\end{center}
\caption{\small Qualitatively correct
evolution of the metric component $g_{2\,2}=m(\rho)$ of
a 3-manifold of revolution under Ricci flow.
The curves have been plotted at unequal time intervals.}
\label{M_PLOT}
\end{figure}

\begin{figure}[tb]
\begin{center}
\includegraphics[width=0.65\textwidth]{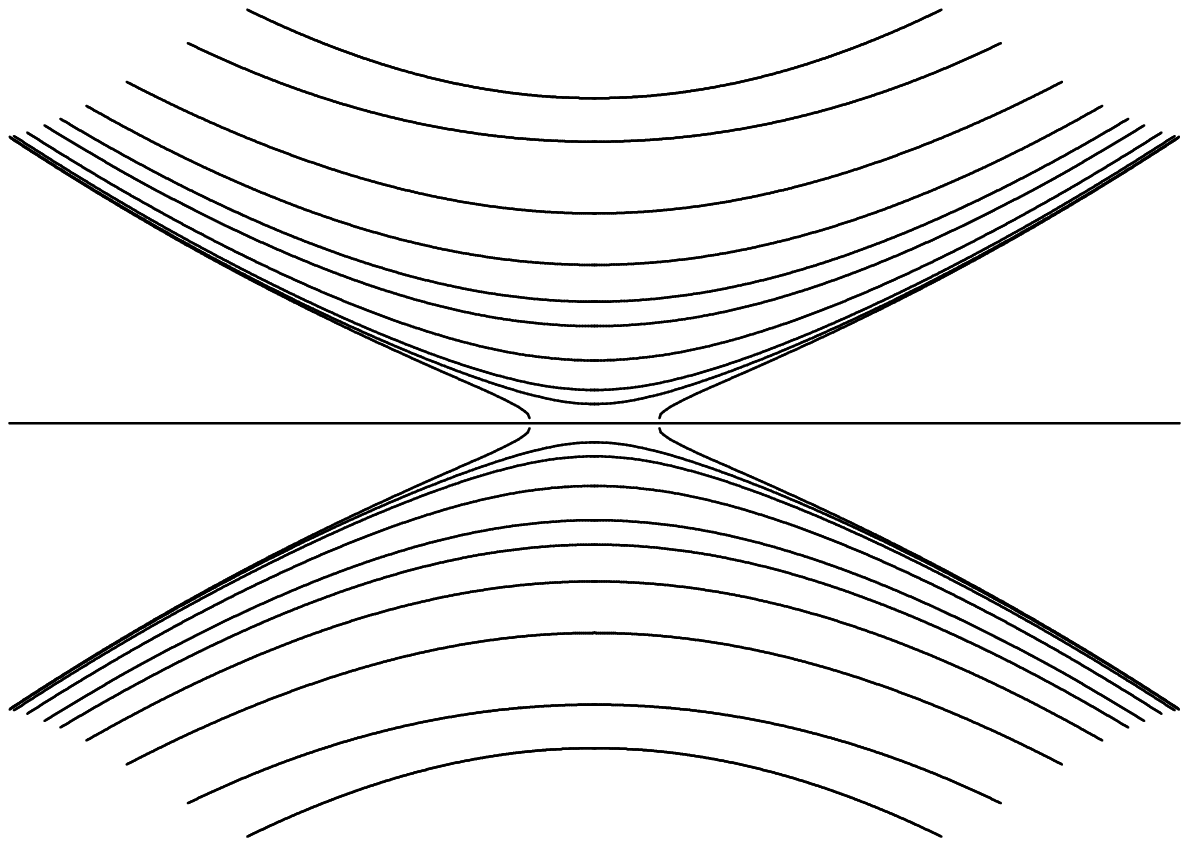}
\end{center}
\caption{\small Qualitatively correct
deformation of the cross-section of
a 3-manifold of revolution under Ricci flow.
The axis of rotation is horizontal.
The pictures have been taken at unequal time intervals.}
\label{CS_PLOT}
\end{figure}

\subsection{A Line cross a Surface of Revolution of Constant Curvature}
\label{HASS}

Let $K_2$ be the Gaussian
curvature of a general abstract Riemannian surface of revolution,
where $K_2$ is any real constant.

The 3-manifold metric is
\begin{equation}
\left[g_{\mu\nu}\right]=\left[\begin{array}{ccc}
\displaystyle h(\rho) & 0 & 0 \\[3ex]
0 & \displaystyle m(\rho) & 0 \\[1ex]
0 & 0 & \displaystyle m(\rho)\,
\cos^2\left(\sqrt{K_2}\,\theta\right)
\end{array}\right],
\end{equation}
where $\rho\equiv x^1$ plays the role of a latitude and
$\theta\equiv x^2$
the role of a latitude on
the abstract Riemannian surface of revolution.

The non-zero elements of the Ricci tensor are
\begin{eqnarray}
R_{1\,1}&=&
\frac{\left(m'\right)^2}{2\,m^2}-\frac{m''}{m}+\frac{m'\,h'}{2\,m\,h} \\[2ex]
R_{2\,2}&=&
\frac{m'\,h'}{4\,h^2}-\frac{m''}{2\,h}+K_2 \\[2ex]
R_{3\,3}&=&
\left(\displaystyle\frac{m'\,h'}{4\,h^2}-\frac{m''}{2\,h}+K_2\right)\,
\cos^2\left(\sqrt{K_2}\,\theta\right),
\end{eqnarray}
and therefore unnormalized Ricci flow satisfies
\begin{eqnarray}
\frac{\partial\,h}{\partial t}&=&
\frac{2\,m''}{m}-\frac{\left(m'\right)^2}{m^2}-\frac{m'\,h'}{m\,h}
\label{URICCI1} \\[2ex]
\frac{\partial\,m}{\partial t}&=&
\frac{m''}{h}-\frac{m'\,h'}{2\,h^2}-2\,K_2. \label{URICCI2}
\end{eqnarray}

The scalar curvature is
\begin{equation}
R=\frac{\left(m'\right)^2}{2\,m^2\,h}-\frac{2\,m''}{m\,h}+\frac{m'\,h'}{m\,h^2}
+\frac{2\,K_2}{m}.
\end{equation}

With respect to the mutually orthogonal unit vectors
\begin{equation}
a=\left(\frac{1}{\sqrt{h}},\, 0,\, 0
\right),\quad
b=\left(
0,\, \frac{1}{\sqrt{m}},\, 0
\right),\quad
\mbox{and}\quad c=\left(
0,\, 0,\, \frac{1}{\sqrt{m}\,\cos\sqrt{K_2}\theta}
\right),
\end{equation}
the sectional curvatures are
\begin{equation}
K(a,b)=K(a,c)=
\frac{\left(m'\right)^2}{4\,m^2\,h}-\frac{m''}{2\,m\,h}+\frac{m'\,h'}{4\,m\,h^2}
\end{equation}
and
\begin{equation}
K(b,c)=\frac{K_2}{m}-\frac{\left(m'\right)^2}{4\,m^2\,h}.
\end{equation}

\subsection{An Example}

We will study the pinching behaviour under unnormalized
Ricci flow of the initial 3-manifold with metric given by
\[
m(\rho)=\frac{1}{10000}+\sin^2\left(\frac{9\,\pi\,\rho}{40}\right)\quad
h(\rho)=1\quad\mbox{and}\quad K_2=1
\]
at $t=0$. We denote the time at which pinching occurs ($m(0)=0$) by
$t=T$.

\subsubsection{Qualitative Behaviour}

The most direct approach to simulating the flow determined by Equations
\ref{URICCI1} and \ref{URICCI2} is to use an explicit finite-difference
formulation. Doing this, we found that the system is (not surprisingly)
highly unstable. We identified two modes of instability. One was of a
short wavelength -- of the type $k\times(-1)^i$ where $i$ is the spatial
grid index and $k$ some constant. This type of instability could be
removed by simply performing the simulation in Maple \cite{HECK}, using
several hundred digits' precision. The other type of instability showed
itself as a long-wavelength smooth perturbation which diverged as the number
of time-steps was increased. We were not able to eliminate this second source
of instability. Practically, it meant that we were restricted
to fairly large
time-steps. To be more precise, we stepped between
$t=0.000025$, $0.000050$, $0.0000625$, $0.00006875$, $0.000071875$,
$0.000075000$,
$0.00007578125$, $0.00007656250$, $0.000076953125$ (shown only
in Figure
\ref{CS_PLOT}) and $0.000077343750$. Since we are therefore not able
to claim that
the method converges to the exact solution, we have chosen to
speak of the
data from this simulation as being only qualitatively correct.

It is nonetheless informative to consider Figures
\ref{H_PLOT}, \ref{M_PLOT} and \ref{CS_PLOT}.

Figures \ref{M_PLOT} and \ref{CS_PLOT} would suggest that pinching
occurs {\it before}
$t=0.000077343750$. At $t=0.000077343750$, values of $m$ at nodes
near $\rho=0$ are
negative. Since Equations \ref{GEN_CURVE} no longer make sense for
negative $m$, an
interval around $\rho=0$ has essentially been cut out. This would of
course correspond
to the view that pinching creates two independent geometric bodies,
as
suggested in Figure \ref{CS_PLOT}. Our single
coordinate $\rho$ should then be replaced by two new coordinates, one
for each body,
and values of $\rho$ which can no longer be attributed to the one or
the other
body should then be discarded as being ``non-geometric''.
It is indeed interesting to consider what might happen {\it after}
pinching, retaining
the coordinate $\rho$ and dropping the condition that the metric
should be
positive definite. It is however important to note that the
negative
values of
$m$ seen in our qualitative numerical results
cannot be trusted, since they are an artifact of Euler's method being
applied with a step-size so large that a singularity was passed over.
What we do claim is
that the cross-sections plotted in Figure \ref{CS_PLOT} do still have
some meaning, if
one does not attempt to read anything into the missing points at the
caps of the two
bodies which resulted from pinching.

We return to this question in Section \ref{NEGATIVE} from another
point of view.

\subsubsection{Asymptotic Scaling}

\begin{figure}[t]
\begin{center}
\includegraphics[width=0.5\textwidth]{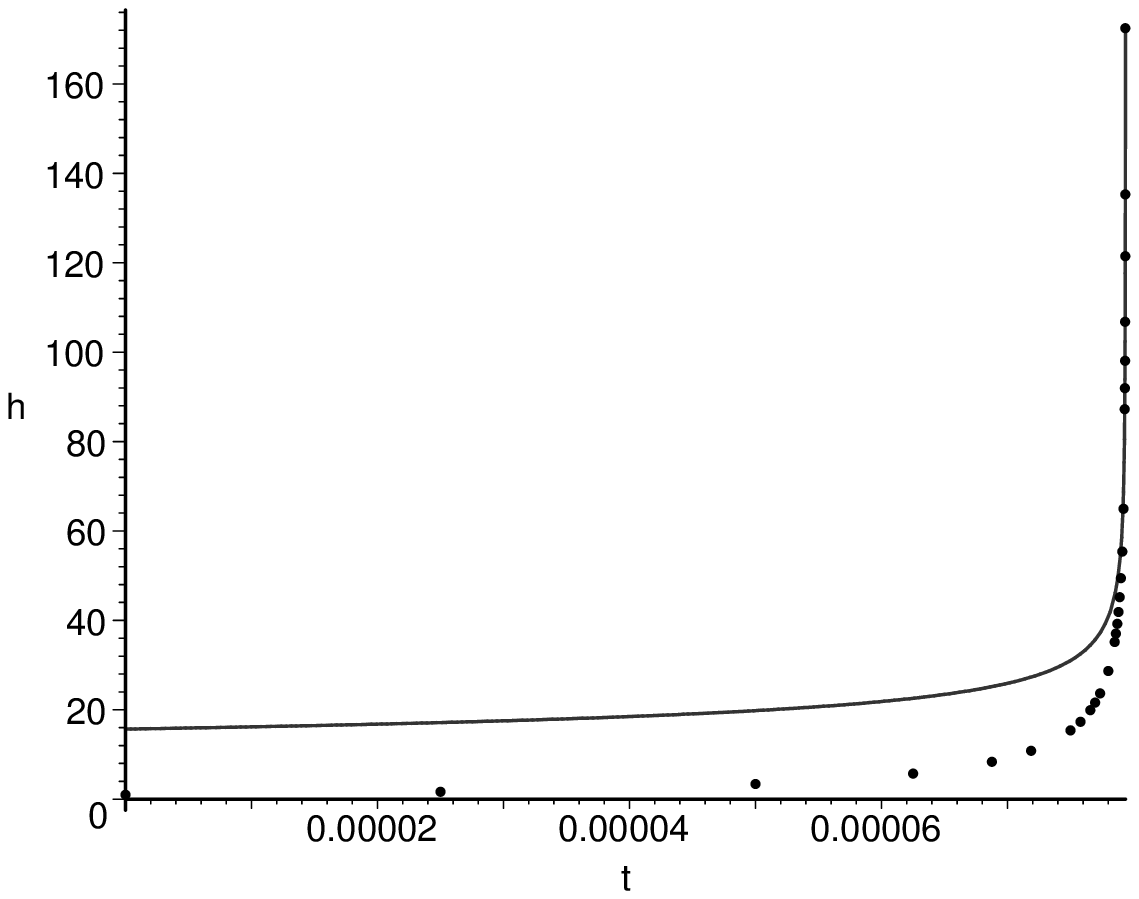}
\end{center}
\caption{\small The evolution of $h$ with time for a $3$-manifold
of revolution.
The points are numerical data. The curve is a fit to the asymptotic
behaviour
of $h$ near pinching.}
\label{ONE}
\end{figure}

\begin{figure}[t]
\begin{center}
\includegraphics[width=0.5\textwidth]{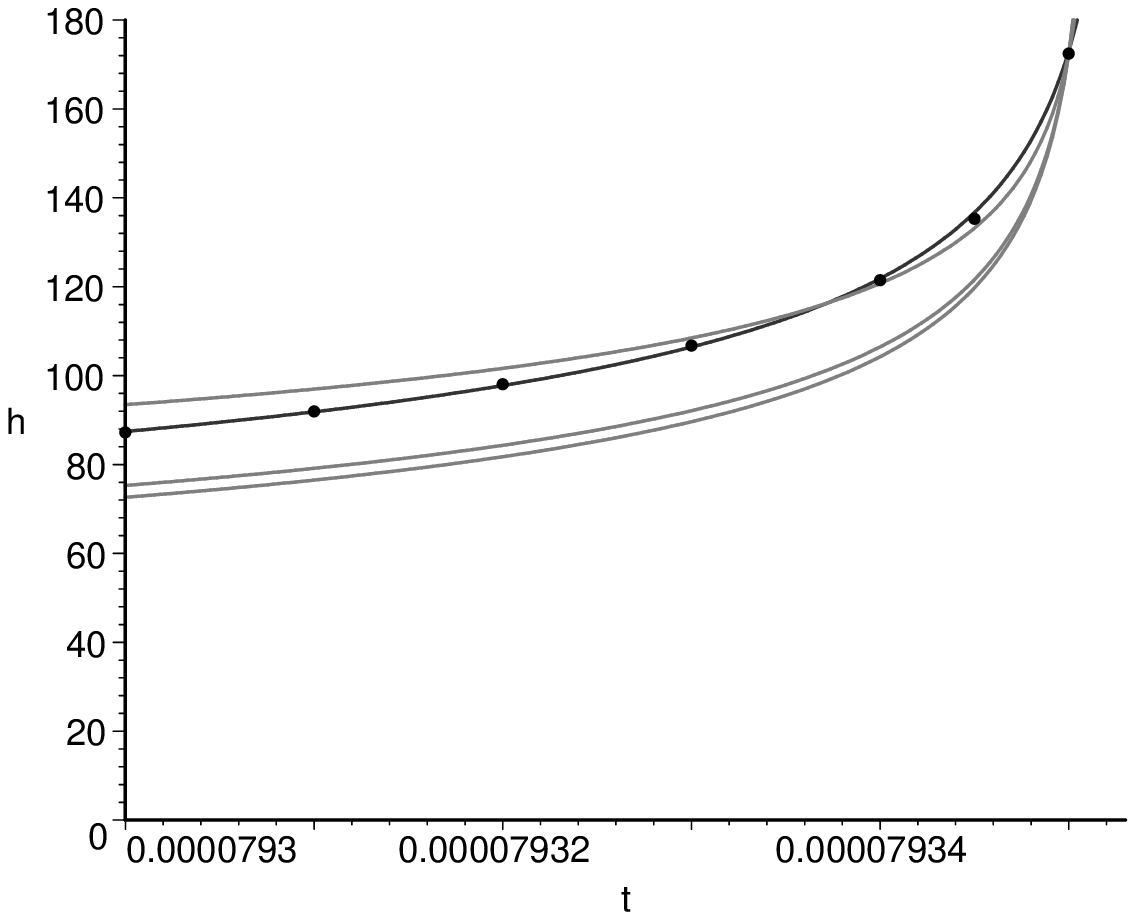}
\end{center}
\caption{\small The evolution of $h$ with time for a $3$-manifold
of revolution close to pinching. The points are numerical data.
The dark curve is the fit
$1.705\times(0.0000793529\!-\!t){}^{-0.235}$. The lighter curves
correspond to other attempts at curve-fitting, in particular
in an attempt to match the pinching time of $T=0.0000793514$.
These are
$1.294\times(0.0000793514\!-\!t)^{-0.24}$,
$1.586\times(0.0000793514\!-\!t)^{-0.23}$, and
$5.389\times(0.0000793514\!-\!t)^{-0.17}$.}
\label{VARI}
\end{figure}

A second approach to studying pinching involves series expansions
of the
metric components $h$ and $m$. We have expanded them to tenth order
in $\rho$
(even terms only) and first
order in $t$:
\[
h(\rho,t)=\sum\limits_{i=0}^{10}(h_i+\dot{h}_i\,t)\,\rho^i
\quad\mbox{and}\quad
m(\rho,t)=\sum\limits_{i=0}^{10}(m_i+\dot{m}_i\,t)\,\rho^i.
\]
Substituting into Equations (\ref{URICCI1}) and (\ref{URICCI2}) we get,
by equating
coefficients (and assuming higher-order terms are zero) equations for
the quantities
$\dot{h}_i$ and $\dot{m}_i$. These equations are quite cumbersome.
We have used
Maple \cite{HECK} both to derive the equations and also to generate C
code \cite{KERNI}
to allow us to compute the flow as rapidly as possible. The C code
was modified to
utilize the {\tt long double} numerical data type, and the code was
run on
a DEC alpha
(the name of a RISC microprocessor) machine
with true 64 bit floating-point arithmetic (the {\tt long double}
type). We
used Euler's method with
a very primitive form of adaptive step-size estimation, finding these
to be satisfactory
for our purposes.

Our data are consistent with the qualitative results of the previous section.

We are able to use these series expansions to make some conjectures
concerning the behaviour of various quantities as $t$ approaches $T$.
The following scaling laws are purely empirical. They are the result
of curve-fitting
to the data we have, in particular for the seven points, at $t=0.000079300$,
$0.000079310$, $0.000079320$, $0.000079330$, $0.000079340$,
$0.000079345$ and $0.000079350$. The step-sizes used were significantly
shorter than the
differences between the times of these representative points.
Figures \ref{ONE} and \ref{VARI} illustrate this fitting process.
It is clear from these illustrations that what we believe to be
asymptotic behaviour sets in fairly late. One may of course ask if
we have fitted to data which is late enough in the flow to truly
warrant being called asymptotic. At this stage, we can only present
the data we have.

All of the following apply to $\rho=0$:
\begin{equation}
\begin{array}{ccr@{}r}
m&\approx&1.409\times(0.0000793514-t)&{}^{0.985} \\[1ex]
h&\approx&1.705\times(0.0000793529-t)&{}^{-0.235} \\[1ex]
R&\approx&0.570\times(0.0000793515-t)&{}^{-1.025} \\[1ex]
K(a,b)&\approx&-1.142\times(0.0000793513-t)&{}^{-0.826} \\[1ex]
K(b,c)&\approx&0.698\times(0.0000793514-t)&{}^{-0.986}
\end{array}
\label{ASYMP}
\end{equation}
It would therefore appear that $T\approx0.0000793514$.

Note that this curve-fitting problem is ill-conditioned,
although Figure \ref{VARI} does suggest that our fits are
as accurate as one could hope for. However,
$h$ appears to
diverge at a significantly later time than our estimate for $T$, and
we do not find
that our curves of best fit for $m$ and $K(b,c)$ obey
$m\times K(b,c)=K_2$
exactly at $\rho=0$.

We do find that the second derivative of the cross-section $y(x)$
(see Equations
\ref{GEN_CURVE}) at $x=0$
diverges as $t\rightarrow T$, indicating that the neck becomes sharper rather
than longer as the surface approaches pinching. This makes sense,
since at $t=T$
we expect to have two abutting caps, as illustrated in Figure \ref{CS_PLOT}.

These results need to be reproduced independently (perhaps via the
DeTurck flow \cite{ISEN}) to ascertain their
validity. We present them on an ``as is'' basis, freely admitting
that we
are not even in a position to give an error analysis,
but in the hope that even
such preliminary data may inspire or support some new analytical
attack.

It is perhaps appropriate to quote Kenk{\=o} at this point
(from essay 82 of \cite{KENKO}):
\begin{quote}
Leaving something incomplete makes it interesting, and gives one the
feeling that there is room for growth. Someone once told me,
``Even when building the imperial palace, they always leave one
place unfinished.''
\end{quote}

One is tempted to think of level set methods (see \cite{SETH} for a
general
introduction) as a possible computational alternative to what we have
tried.
Level set methods have been successfully applied to mean-curvature flow
\cite{MEAN}, and indeed their strength lies in their ability to cope
naturally with changes in topology of an evolving surface, as
pinching involves.

We feel that some reformulation of the problem, perhaps along the
lines of what follows, may facilitate a level set formulation of Ricci
flow, and that this is without doubt a worthwhile focus for further work.

\section{Speculation}
\label{NEGATIVE}

The system of differential equations (\ref{URICCI1}) and (\ref{URICCI2})
above makes sense when the
metric tensor is not necessarily positive definite, so long as
smoothness occurs
at points where $m=0$. So we are lead to speculate about the idea
of Ricci flow
in the larger space of positive definite {\it and} indefinite metrics.

The advantage of such an approach might be that instead of the
manifold pinching, one might allow that part of the
manifold has a positive definite metric and part has an indefinite
metric, after a singularity occurs.
In this case, the computational problem might become more tractible.
The reason is that the flowing family of indefinite metrics
would always have a fixed domain of definition.

In particular, at a neck pinch, the initially Riemannian metric for a
rotationally symmetric sphere
should ``jump'' to a metric with a positive definite part
on two polar caps and an indefinite part on an equatorial band. The
indefinite part will expand and
the positive definite parts contract until the metric on the whole space
becomes
indefinite.
Further numerical experimentation might
give a better idea of a suitable
formulation for this process.

Note that Perelman \cite{PEREL1} has indicated how his notion of entropy
can be used as a potential function for Ricci flow. Extending this
approach to the larger domain of positive definite and indefinite 
metrics may lead to a
true level set formulation of Ricci flow, in this space of metrics on 
a fixed manifold.
Note that at singularities where some components of the Ricci tensor 
become infinite, one would expect the
level sets to fatten ( see e.g. \cite {ALT}).


\begin{thebibliography}{99}
\bibitem{ALT} S. Altschuler, S. Angenent and Y. Giga,
{\it Mean curvature flow through singularities for surfaces of revolution},
The Journal of Geometric Analysis, Volume 5.  293--358 (1995).
\bibitem{ANG} S. Angenent and D. Knopf, {\it An example of neck
pinching for Ricci Flow on $S^{n+1}$}, to appear in Mathematical
Research Letters (2004).
\bibitem{BH} S. Bleiler and C. Hodgson, {\it Spherical space forms
and Dehn filling}, Topology, Volume 35, 809--833, (1996).
\bibitem{SPECTRAL} G.L. Browning, J.J. Hack and P.N. Swartztrauber,
{\it A Comparison of Three Numerical Methods for Solving Differential
Equations on the Sphere}, Monthly Weather Review, Volume 117,
1058--1075, (1989).
\bibitem{CAO} H. Cao, B. Chow, S. Chu and S.T. Yau, Editors, {\it
Collected Papers on Ricci Flow}, Series in Geometry and Topology,
Volume 37, International Press, MA (2003).
\bibitem{MEAN} D.L. Chopp and J.A. Sethian,
{\it Flow under Curvature: Singularity Formation, Minimal Surfaces,
and Geodesics},
Experimental Mathematics, Volume 2 (1993) 235--255.
\bibitem{BCHOW} B. Chow, {\it The Ricci Flow on the 2-Sphere},
in H. Cao et al., Editors,
{\it Collected Papers on Ricci Flow}, Series in Geometry and Topology,
Volume 37, International Press, MA (2003) 226--237.
\bibitem{CHOW} B. Chow and D. Knopf, {\it The Ricci Flow: An
Introduction}, Mathematical Surveys and Monographs, Volume 110,
AMS (2004).
\bibitem{ENGMAN} M. Engman, {\it A Note on Isometric Embeddings of
Surfaces of Revolution}, American Mathematical Monthly,
Volume 111 (2004) 251--255.
\bibitem{ISEN} D. Garfinkle and J. Isenberg,
{\it Critical behavior in Ricci flow}, arXiv:math.DG/0306129.
\bibitem{GT} M. Gromov and W. Thurston, {\it Pinching constants for
hyperbolic manifolds}, Invent. Math. Volume 89, 1--12, (1989).
\bibitem{HECK} A. Heck, {\it Introduction to Maple}, Third Edition,
Springer-Verlag, New York, (2003).
\bibitem{IVEY} T. Ivey, {\it The Ricci Flow on Radially Symmetric
$R^3$}, Comm. in P.D.E.,
Volume 19, 1481--1500, (1994).
\bibitem{KENKO} D. Keene, Translator,
{\it Essays in Idleness: The Tsurezuregusa of Kenk{\=o}},
Columbia University Press, New York (1998).
\bibitem{KERNI} B.W. Kernighan and D.M. Ritchie,
{\it The C Programming Language}, Second Edition,
Prentice Hall, New Jersey (1988).
\bibitem{HYAM} P.R.A. Leviton and J.H. Rubinstein,
{\it Deforming Riemannian Metrics on the 2-Sphere}, in
L. Simon and N.S. Trudinger, Editors, {\it Miniconference on
Geometry and Partial Differential Equations (Canberra 1985)},
Proceedings of the Centre for Mathematical Analysis ANU,
Volume 10, Australian National University, Canberra (1986) 123--127.
\bibitem{OPEN} OpenGL Architecture Review Board et al.,
{\it OpenGL Programming Guide}, Fourth Edition,
Addison-Wesley, MA (2003).
\bibitem{PEREL1} G. Perelman, {\it The Entropy Formula for the Ricci
Flow and its geometric applications}, arXiv:math.DG/0211159.
\bibitem{PEREL2} G. Perelman, {\it Ricci Flow with surgery on
three-manifolds}, arXiv:math.DG/0303109.
\bibitem{SETH} J.A. Sethian, {\it Level Set Methods}, Cambridge
Monographs on Applied
and Computational Mathematics Vol 3, Cambridge University Press,
Cambridge (1996).
\bibitem{SIMON} M. Simon, {\it A class of Riemannian manifolds that
pinch when evolved by Ricci flow}, Manuscripta Math. Volume 101,
89--114, (2000).
\end{thebibliography}
\end{document}